\documentclass{article}
\usepackage{amssymb,amsthm,amsmath,graphicx}
\usepackage{amsfonts}
\newtheorem{theorem}{Theorem}[section]

\newtheorem{remark}{Remark}[section]

\begin{document}

\title{\bf{Influence Function Analysis of the Restricted Minimum Divergence Estimators : A General Form }}
\author{
       \\  
       Abhik Ghosh\\
       Indian Statistical Institute, Kolkata, India.\\
}
\date{}
\maketitle

\begin{abstract}
The minimum divergence estimators have proved to be 
useful tools in the area of robust inference. The robustness of such estimators are measured
using the classical Influence functions. However, in many complex situations like 
testing a composite hypothesis using divergence require the estimators to be restricted 
into some subspace of the parameter space. The robustness of these restricted 
minimum divergence estimators are very important in order to have overall robust inference. 
In this paper we provide a comprehensive description of the robustness of such restricted estimators
in terms of their Influence Function for a general class of density based divergences 
along with their unrestricted versions. In particular, the robustness
of some popular minimum divergence estimators are also demonstrated under certain usual restrictions. 
Thus this paper provides a general framework for the influence function analysis of a large class of  
minimum divergence estimators with or without restrictions on the parameters.
\end{abstract}


\section{Introduction}
 
The minimum divergence approach has proved to be a very useful one in the context 
of parametric statistical inference. The idea behind this approach is 
to quantify the discrepancy between the sample data and the parametric 
model through an appropriate divergence and minimize this discrepancy measure over 
the parameter space. There are two ways of such quantification in literatures; --
 either through the distribution functions or the probability density functions 
 (with respect to some suitable measure).  Most of the density based minimum divergence methods 
 are seen to be particularly useful due to their strong robustness along with high efficiencies.

However, in many complex statistical problems we need to estimate the parameter of interest 
under some pre-specified restrictions on the parameter space. For example, when testing a 
composite hypothesis, we need to estimate the parameter under the restriction imposed by null hypothesis. 
For such cases we need to minimize the divergence measures only over a restricted subspace of the  parameter space.
Simpson (1989), Lindsay (1994) and Basu et al.~(2013) used such restricted minimum divergence estimators in the 
context of testing statistical hypothesis and derived their asymptotic properties. 
But they did not consider the robustness of these restricted estimators through the usual indicators, 
although it is also very important in order to obtain robust solution for the overall complex inference problem. 
Indeed, the robustness aspect of the restricted minimum divergence estimators are not 
well studied in literatures. In this paper, we will consider this very important issue 
and describe the robustness of the general minimum divergence estimators in terms of the 
Influence Function Analysis.

The rest of this paper is organized as follows. Sections \ref{SEC:MDE_IF} describes the concept of the 
minimum divergence estimators and present a general form for their influence 
function analysis. In Section \ref{SEC:IF_RMDE_general}, we will 
derive a general form of the influence function of the restricted minimum divergence estimators.
Finally, in Section \ref{SEC:part_div}, we will apply the general results in case of some popular divergences - disparity,
density power divergence, S-divergence; under certain usual restrictions on the parameter of interest.

\section{Density-Based Minimum Divergence Estimation and  Influence Function : General Form   }\label{SEC:MDE_IF}

Let us begin our discussion with a general parametric estimation problem. 
We have $n$ independent and identically distributed observations $X_1,\ldots, X_n$ from a 
distribution $G$. We want to model it by a parametric family of distributions  
${\cal F}_\theta =\{F_\theta: \theta \in \Theta \in {\mathbb R}^p\}$. 
Without loss of generality, let the support of $G$ 
and the parametric model ${\cal F}_\theta$ are the same. 
Also let both $G$ and ${\mathcal F}_\theta$ 
belong to $\mathcal{G}$, the class of all distributions having densities 
with respect to the appropriate $\sigma-$finite measure $\mu$ on
 the $\sigma-$field $(\Omega, \mathcal{A})$ and $f_\theta$, $g$ be the 
density functions of $F_\theta$, $G$ respectively with respect to $\mu$. 
We want to estimate the parameter $\theta$ based on the available sample data.
In case of density-based minimum divergence estimation, this is done by choosing the model 
element that provides the closest match to the data where the separation 
between the model and data is quantified by a nonnegative function $\rho(\cdot,\cdot)$ from
 $\mathcal{G} \times \mathcal{G}$ to $[0,~\infty)$ that equals zero if and only if its arguments 
 are identically equal. Such functions $\rho(\cdot,\cdot)$ are termed as the Statistical Divergence 
 and the estimator $\hat{\theta}$  of $\theta$ obtained by minimizing $\rho(\hat{g},f_\theta)$ 
 with respect to $\theta \in \Theta$, where $\hat{g}$ is some nonparametric estimator
  of $g$ based on the sample data, is called the Minimum Divergence Estimator (MDE).
    In terms of statistical functionals, the Minimum Divergence Functional $T_\rho(G)$ corresponding
   to the divergence $\rho(\cdot,\cdot)$ is defined  by the relation
$\rho(g, f_{T_\rho(G)}) = \min_{\theta \in \Theta} ~ \rho( g, f_\theta)$ 
provided such a minimum exists.
There are several popular examples of the statistical divergences that generate highly robust 
  and efficient estimators including the Disparity family (Lindsay, 1994), Cressie-Read Power Divergences 
  family (Cressie and Read, 1984), Density Power Divergence (Basu et al., 1998) etc. See 
  Csiz\'{a}r (1963, 1967a, b), Ali and Silvey (1966), Vajda (1989), Pardo (2006) and
   Basu et al.~(2011) for further examples and the details of the 
  minimum divergence estimators including there asymptotic and robustness properties. 
  Most of the divergences used for statistical inference have particular form, given by 
  \begin{equation}
\rho( g, f) = \int D(g,f) d\mu 
\label{EQ:part_div}
  \end{equation}
for some suitable function $D(\cdot,\cdot) : \mathbb{R}\times\mathbb{R} \longmapsto [0~\infty)$. 
So, in this paper also, we will restrict our attention to the divergences satisfying Equation 
(\ref{EQ:part_div}) only. Then, the estimating equation of the MDE is given by
 \begin{equation}\label{EQ:estimating_equation_mde}
 \nabla  \rho(g, f_{\theta}) = \int \nabla D(g,f_{\theta}) d\mu  = 0,
\end{equation} 
where $\nabla$ represents the derivative with respect to $\theta$. Note that, this estimating equation 
does not necessarily give us an M-estimator; it does so only when $\nabla D(g,f_{\theta})$
containing $f_\theta$ includes only the linear function of $g$ or some constant independent of $g$. 
However, the number of divergences satisfying this condition is limited (See, eg.~Patra et al., 2013) 
so that we can not always apply the theory of M-estimators to describe the properties of the MDEs. 
However, all the MDEs obtained as a solution to (\ref{EQ:estimating_equation_mde}) will be Fisher consistent 
by definition of $\rho(\cdot,\cdot)$.

	The MDEs are mostly popular due to their strong robustness and 
	in this context a useful tool is the Influence Function (Hampel,1968, 1974)which is an indicator of their
	classical first-order robustness, as well as of their asymptotic efficiency. 
	To obtain the influence function of the minimum divergence estimators based on 
	the divergence $\rho(\cdot,\cdot)$, we consider the $\epsilon$ contaminated version 
	of the true density $g$ given by $g_\epsilon(x) = (1- \epsilon)g(x) + \epsilon \chi_y(x)$. 
	Similarly $G_\epsilon(x) = (1-\epsilon)G(x) + \epsilon \wedge_y(x)$; here $\chi_y(x)$ and $\wedge_y(x)$ 
	are respectively density and distribution function of the degenerate distribution at $y$.
	 Let $\theta^g = T_\rho(G)$  and $\theta_\epsilon = T_\rho(G_\epsilon)$ be the functional obtained via the minimization 
	 of $ \rho(g, f_\theta)$ and $ \rho(g_\epsilon, f_\theta)$ respectively. 
	 Then the Influence function of the Minimum Divergence Functional $T_\rho(\cdot)$ 
	 is defined as
$IF(y, T_\rho, G) = \frac{\partial \theta_\epsilon}{\partial \epsilon}\big|_{\epsilon=0}.$
But from the definition of $\theta_\epsilon$, 
it must satisfy the estimating equation (\ref{EQ:estimating_equation_mde}).
Now substituting $g_\epsilon$ and $\theta_\epsilon$ in place of $g$ and $\theta$ 
in (\ref{EQ:estimating_equation_mde}) respectively and 
differentiating with respect to $\epsilon$ at $\epsilon=0$ we get ,
\begin{eqnarray}
&&
\int \frac{\partial[\nabla D(g(x),f_{\theta^g}(x))] }{\partial g(x)}[- g(x) + \chi_y(x)] d\mu(x) 
\nonumber \\ && ~~~~~
+ \int \frac{\partial[\nabla D(g(x),f_{\theta^g}(x))] }{\partial f_{\theta^g}(x)} [\nabla f_{\theta^g}(x)]^T 
IF(y, T_\rho, G)  d\mu(x) = 0. \nonumber 
\end{eqnarray}
But $$\nabla D(g(x),f_\theta(x)) = D^{(2)}(g(x),f_\theta(x)) \nabla f_\theta(x);$$ 
hence
$$ \frac{\partial[\nabla D(g(x),f_\theta(x))] }{\partial g(x)} = D^{(1,2)}(g(x),f_\theta(x)) \nabla f_\theta(x)$$
and
$$\frac{\partial[\nabla D(g(x),f_\theta(x))] }{\partial f_\theta(x)} 
= D^{(2)}(g(x),f_\theta(x))  \frac{\partial[\nabla f_\theta(x)] }{\partial f_\theta(x)}  
+ D^{(2,2)}(g(x),f_\theta(x)) \nabla f_\theta(x).$$ 
Here, $D^{(i)}(\cdot,\cdot)$ denotes the first order partial derivative of $D(\cdot,\cdot)$ 
with respect to its $i^{th}$ argument, $D^{(i,j)}(\cdot,\cdot)$ denotes its second order partial derivative 
with respect to  $i^{th}$  and $j^{th}$ arguments ($i, j =1,2$) and we have assumed that the standard regularity 
conditions hold for the densities so that all above derivatives exists and can be interchanged with the integrals. 
Thus the expression of the influence function of the minimum divergence functional $T_\rho$ simplifies to,
\begin{eqnarray}\label{EQ:IF_MDE}
 IF(y, T_\rho, G)  = N(\theta^g)^{-1} [\xi(\theta^g) - M(y;\theta^g)],
\end{eqnarray}
%
where 
$$N(\theta) = \int \left[D^{(2)}(g(x),f_\theta(x)) \nabla_2 f_\theta(x) + D^{(2,2)}(g(x),f_\theta(x)) \{\nabla f_\theta(x)\} \{\nabla f_\theta(x)\}^T\right]  d\mu(x),$$
$$M(y;\theta)= D^{(1,2)}(g(y),f_\theta(y)) [\nabla f_\theta(y)],$$
and
$$\xi(\theta)= \int  D^{(1,2)}(g(x),f_\theta(x)) [\nabla f_\theta(x)] g(x) d\mu(x) 
= E_g\left[M(X;\theta)\right].$$

In particular, when the true distribution
$G$ belongs to the parametric model, so that the density $g(x) = f_{\theta_0}(x)$ 
for some $\theta_0 \in \Theta$, we get $\theta^g = \theta_0$ and the influence function becomes 
$IF(y, T_\rho, F_{\theta_0})  = N(\theta_0)^{-1} [\xi(\theta_0) - M(y;\theta_0)]$. Therefore, 
the influence functions of the MDEs will be bounded at the model for all those divergences 
for which the function $\left|M(y;\theta)\right|$ is bounded in $y$ for all $\theta$.

	Further, note that as expected from the interpretation of the influence function by Hampel et al.~(1986),
	we have 
$$\int IF(y, T_\rho, G)  d(G(y)) = \int IF(y, T_\rho, G)  g(y) d\mu(y) = 0.$$ 
Also, if the MDE $T_n = T(\hat{G})$ is an $\sqrt n$-estimator of $T(G)$, then it follows that (Hampel et al.~)
the asymptotic distribution of $\sqrt{n} (T_n - T(G))$ is asymptotically normal with mean zero and variance 
$$V(T,G) = 	\int IF(y, T_\rho, G)IF(y, T_\rho, G)^T  d(G(y)) 
=  N(\theta^g)^{-1} Var_g[M(X;\theta^g)] N(\theta^g)^{-1},$$
where $Var_g(\cdot)$ denotes the variance under the distribution of $g$.

\section{The Influence Function of Restricted MDE : General Case }\label{SEC:IF_RMDE_general}

We will now consider the case of restricted minimum divergence estimators and 
derive a general expression for its influence function extending the concepts of the previous section. 
Consider the set-up of the previous Section \ref{SEC:MDE_IF}, but now we want to estimate
the parameter $\theta$ only  over a restricted (proper) subspace $\Theta_0$ of the whole
 parameter space $\Theta$. In most of the cases, we can define the subspace $\Theta_0$ by
 a set of $r$ restrictions of the form 
 \begin{equation}\label{EQ:restriction}
 h(\theta)=0 ~~~~~~~~ \mbox{on } ~ \Theta,
 \end{equation}
 for some function $h : {\mathbb R}^p \longmapsto {\mathbb R}^{r}$ satisfying the property that
 the $p\times r$ matrix 
$
H(\theta)= \frac{\partial h(\theta)}{\partial \theta}
$ 
exists with rank $r$ and is continuous in $\theta$. Thus, under $\Theta_0$, the parameter 
 $\theta$ essentially contains $p - r$ independent parameters.

	We can solve the above estimation problem 
 by minimizing $\rho(\hat{g},f_\theta)$  with respect to $\theta \in \Theta_0$ and the estimator
  obtained from this minimization exercise will be called the Restricted Minimum Divergence 
  Estimator (RMDE). Lindsay (1994) and Basu et al.~(2013) derived the asymptotic distribution of 
  such RMDE for the disparity family and the density power divergences respectively.
 Let us define the Restricted Minimum Divergence 
 Functional $\widetilde{T}_\rho(G)$ by the relation
 $
 \rho(g, f_{\widetilde{T}_\rho(G)}) = \min_{\theta \in \Theta_0} ~ \rho( g, f_\theta) 
 = \min_{h(\theta) = 0} ~ \rho( g, f_\theta),
 $
 provided such a minimum exists. We can easily solve this minimization problem using the
  Lagrange multiplier method.
  
Now to derive the influence function of the Restricted Minimum divergence functional, as before,
 we will consider the $\epsilon$-contaminated density $g_\epsilon(x)$ and let 
 $\widetilde{\theta^g }= \widetilde{T}_\rho(G)$  and 
 $\widetilde{\theta_\epsilon } = \widetilde{T}_\rho(G_\epsilon)$. 
 Note that $\widetilde{\theta_\epsilon }$ is the minimizer of $ \rho(g_\epsilon, f_\theta)$ subject to 
 (\ref{EQ:restriction}). Let us consider the restrictions which can be substituted explicitly in the 
 expression of $ \rho(g_\epsilon, f_\theta)$ before taking its derivatives with respect to $\theta$;
 the corresponding derivative will be then zero at $\theta=\widetilde{\theta_\epsilon }$ and 
 proceeding as in Section \ref{SEC:MDE_IF}, we get 
 \begin{equation}
N_0(\widetilde{\theta^g }) IF(y, \widetilde{T}_\rho, G) - \xi_0(\widetilde{\theta^g }) 
+ M_0(y;\widetilde{\theta^g }) = 0, \label{EQ:IF_eq1}
 \end{equation}
where $N_0(\theta)$, $\xi_0(\theta)$, $M_0(y;\theta)$ are the same as $N(\theta)$, $\xi(\theta)$, $M(y;\theta)$ 
respectively but with an additional restriction of $h(\theta)=0$. 
Also, since $\widetilde{\theta_\epsilon }$ must satisfy (\ref{EQ:restriction}), 
a differentiation with respect to $\epsilon$ at $\epsilon=0$ yields
\begin{equation}
H(\widetilde{\theta^g })^T IF(y, \widetilde{T}_\rho, G) = 0, \label{EQ:IF_eq2}
\end{equation}
We need to solve the two equations (\ref{EQ:IF_eq1}) and (\ref{EQ:IF_eq2}) 
to get a general expression for the influence function  $IF(y, \widetilde{T}_\rho, G)$. 
Combining them, we get
  \begin{equation}
 \begin{pmatrix}
 {N_0}(\widetilde{\theta^g }) \\ 
 H(\widetilde{\theta^g })^T
 \end{pmatrix} IF(y, \widetilde{T}_\rho, G) = \begin{pmatrix}
  \xi_0(\widetilde{\theta^g }) - M_0(y;\widetilde{\theta^g }) \\ 
  {\mathbf 0}
  \end{pmatrix}. 
    \end{equation}
After simplification, we get the general expression for the influence function of 
Restricted Minimum Divergence functional which is presented in the following Theorem:

  \begin{theorem}\label{THM:IF_RMDE}
Consider above notations and assume that $rank(H(\widetilde{\theta^g })) = r $. 
Then the Influence Function of the Restricted Minimum Divergence Equation corresponding to 
the divergence (\ref{EQ:part_div}) is given by
  \begin{equation}
IF(y, \widetilde{T}_\rho, G) = \left[{N_0}(\widetilde{\theta^g })^T {N_0}(\widetilde{\theta^g }) + 
   H(\widetilde{\theta^g }) H(\widetilde{\theta^g })^T\right]^{-1} {N_0}(\widetilde{\theta^g })^T 
   \left[ \xi_0(\widetilde{\theta^g }) - M_0(y;\widetilde{\theta^g }) \right],
  \end{equation}
  provided ${N}_0(\cdot)$, $\xi_0(\cdot)$ and $M_0(\cdot)$ can be defined as above.
  \qed
  \end{theorem}

In particular, if the true density belongs to the model family and the imposed restrictions are valid, i.e., 
$g=f_{\theta_0}$ for some $\theta_0$ satisfying $h(\theta_0)=0$, then we just put 
$\widetilde{\theta^g }=\theta_0$ to obtain the corresponding influence function.
Therefore, the influence functions of the RMDEs will be bounded at the model for all those divergences 
for which the function $\left|M_0(y;\theta)\right|$ is bounded in $y$ for all $\theta$. in particular, whenever
the IF of the MDE at the model is bounded the IF of RMDE at the model will also be bounded at the model for any restrictions; but the converse is not true.

\begin{remark}\label{REM:asymp_RMDE}
It is easy to check that  
$$
\int IF(y, \widetilde{T}_\rho, G)  d(G(y)) = \int IF(y, \widetilde{T}_\rho, G)  g(y) d\mu(y) = 0.
$$ 
Thus, if the RMDE $\widetilde{T}_n = \widetilde{T}(\hat{G})$ is an $\sqrt n$-estimator of $\widetilde{T}(G)$, 
then it follows that (Hampel et al., 1986)
the asymptotic distribution of $\sqrt{n} (\widetilde{T}_n - \widetilde{T}(G))$ is asymptotically normal with mean zero and variance 
$$V(\widetilde{T},G) = 	\int IF(y, \widetilde{T}_\rho, G)IF(y,\widetilde{T}T_\rho, G)^T  d(G(y)).$$
At the model $g=f_{\theta_0}$ for some $\theta_0$ satisfying $h(\theta_0)=0$, above expression of asymptotic variance further simplifies to 
\qed
\end{remark}

We will now explore a couple of particular cases of restrictions that are commonly used in parametric estimation.
\\

\textbf{Example 1:}
First we will consider a simple and perhaps most popular case of restrictions where few components of the parameter $\theta$ is pre-specified.
Precisely, let $\theta=(\theta_1 ~ \theta_2)^T$ where $\theta_1$ is an $r$-vector and its value is specified at
 $\theta_{1,0}$ as restrictions. Thus we consider the RMDE of $\theta$ under the restriction
  $\theta_1 = \theta_{1,0}$. Note that, intuitively, in this case we must have RMDE of $\theta_1$ to be fixed at $\theta_{1,0}$ having zero influence function and the influence function analysis of the RMDE of $\theta_2$ should be 
  the same as that of the unrestricted MDE considering $\theta_2$ as the only parameter of interest. 
  We will now apply the general formulas derived above to this simple case to verify if those general results are in-line with the intuitive results.
  
 First consider the influence function of the MDE of $\theta$ in the unrestricted case given by 
  $IF(y, T_\rho, G)  = N(\theta^g)^{-1} [\xi(\theta^g) - M(y;\theta^g)]$. Let us partition this result
   in terms of $\theta_1$ and $\theta_2$ to get    $ M(y;\theta) = (M_1(y;\theta) ~ M_2(y;\theta))^T$,
   $\xi(\theta) = (\xi_1(\theta) ~ \xi_2(\theta))^T$ and 
  $$
   N(\theta) = \begin{pmatrix}
    N_{11}(\theta) & N_{12}(\theta) \\ 
    N_{12}(\theta)^T & N_{22}(\theta) 
    \end{pmatrix}
  $$
  where $M_1$ and $\xi_1$ are $r$-vectors and $N_{11}$ is the matrix of order $r \times r$.  
 
 	Now, considering the restricted case, we have $h(\theta)=\theta_1 - \theta_{1,0}$ so that
 	$$
 	H(\theta) = \begin{pmatrix}
 	    I_r \\ 
 	    O 
 	    \end{pmatrix}.
 	$$
Also, $\widetilde{\theta^g } = (\theta_{1,0}, ~ \theta_2^g)^T$ and hence
 	$$
 	N_0(\widetilde{\theta^g }) =\begin{pmatrix}
 	   O_{r\times r} & O \\ 
 	    O & N_{22}((\theta_{1,0}, ~ \theta_2^g)^T) 
  		\end{pmatrix},
  	$$
   $ M_0(y;\widetilde{\theta^g }) = [0_r ~ M_2(y;(\theta_{1,0}, ~ \theta_2^g)^T)]^T$ and 
     $\xi_0(\widetilde{\theta^g }) = [0_r ~ \xi_2((\theta_{1,0}, ~ \theta_2^g)^T)]^T$.Therefore , using the above result, the Influence function of the RMDE of $\theta$ becomes
     $$
     IF(y, \widetilde{T}_\rho, G)  = \begin{pmatrix}
     0_r \\
     N_{22}((\theta_{1,0} \theta_2^g)^T)^{-1} \left[\xi_2((\theta_{1,0} \theta_2^g)^T) -  M_2(y;(\theta_{1,0} \theta_2^g)^T)\right].
     \end{pmatrix}
     $$
     Thus the influence function of the RMDE corresponding to $\theta_1$ is zero and that corresponding to $\theta_2$ 
     is the same as that obtained in the unrestricted case considering $\theta_2$ only, as expected.
\hfill{$\square$}

\begin{remark}
Note that Above Theorem \ref{THM:IF_RMDE} can only be applied provided the restrictions are such that 
$rank(H(\widetilde{\theta^g })) = r $. But in many practical situations we need to consider restrictions for which 
the rank is strictly less than $r$ and we can not apply the above Theorem \ref{THM:IF_RMDE} directly to obtain 
the influence function of the corresponding RMDEs. However, the  arguments presented to derive the theorem can
 still be applied with some small modifications as required. One such common case is presented below in Example 2.
\hfill{$\square$}
\end{remark}

\textbf{Example 2:}
Let us now consider another slightly complicated case of restrictions where the first $r$ components 
of $\theta$ depend among themselves through only one unknown parameter, say $\beta$. 
Such restrictions are common in case of multivariate normal models with mean $\mu$ and variance $\sigma^2 I_p$
 when we consider the restrictions $\mu = \beta \mu_0$ with known $\mu_0$. 
 And estimation of the parameter $\beta$ and $\sigma^2$ are important under such restrictions 
 for various composite testing problems with $p$ independent normal populations. 
For example, while testing for homogeneity of mean among the $p$ normal populations with unknown equal variances, 
we have to consider the specified restrictions with $\mu_0=(1, \cdots, 1)^T$ under the null hypothesis. 
 In general, let $\theta=(\theta_1 ~ \theta_2)^T$
  where $\theta_1$ is an $r$-vector and assume that $\theta_1 = \phi(\beta)$ with known function 
  $\phi : \mathbb{R} \mapsto \mathbb{R}^{r}$ . We will assume that $\phi(\beta) = (\phi_1(\beta), \cdots, \phi_r(\beta))^T$ and each $\phi_i$ are twice differentiable real functions with non-zero 
  derivatives. Here also we will consider the partitions of the matrices  $N(\theta)$, $\xi(\theta)$
  and $M(y;\theta)$  in terms of $\theta_1$ and $\theta_2$ as in Example 1.
  
  To derive the influence function of the RMDEs in this case, note that $h(\theta) = \theta_1 - \phi(\beta)$ so that 
  $$
   	H(\theta) = \begin{pmatrix}
   	   I_r - B \\ 
   	    O 	\end{pmatrix},
   $$
  where the $r \times r$ matrix $B$ is defined as $B = \frac{\partial \phi(\beta)}{\partial \theta_1}$.
  Note that, the $(i, j)^{th}$ element of the matrix $B$ is given by 
  $b_{ij} = \frac{\phi_j'(\beta)}{\phi_i'(\beta)}$ for each $i, j = 1, \cdots, r$, 
  where $\phi_i'(\beta)$ is the first derivatives with respect to $\beta$. 
Next, simple differentiation gives that, 
$$
\nabla  f_{(\phi(\beta), ~\theta_2)^T}(x) = \begin{pmatrix}
     	   B \frac{\partial f_\theta(x)}{\partial \theta_1} \\ 
     	   \frac{\partial f_\theta(x)}{\partial \theta_2} \end{pmatrix} = B^* \nabla f_\theta(x),
$$  
where, $B^*$ is a $p \times p$ matrix defined as
$$
       	B^* = \begin{pmatrix}
       	    B & O\\ 
       	    O & I	\end{pmatrix},
$$  
and
$$
\nabla_2  f_{(\phi(\beta), ~\theta_2)^T}(x) = \begin{pmatrix}
B \frac{\partial^2 f_\theta(x)}{\partial \theta_1^2}B^T 
+ B^{(1)}\left[\frac{\partial f_\theta(x)}{\partial \theta_1} \otimes I_r\right] &  
B \frac{\partial^2 f_\theta(x)}{\partial \theta_2\partial \theta_1}\\ 
\frac{\partial f_\theta(x)}{\partial \theta_1 \partial \theta_2} B^T & 
\frac{\partial^2 f_\theta(x)}{\partial \theta_2^2} \end{pmatrix},
$$ 
where the $r \times r^2$ matrix $B^{(1)}$ is defined as 
$B^{(1)} = \frac{\partial^2 \phi(\beta)}{\partial \theta_1^2}$.
Then we have $M_0(y;\theta) = B^* M(y;\theta)$, $\xi_0(\theta) = B^* \xi(\theta)$ and
$$
N_0(\theta) = B^* M(\theta)(B^*)^T + \begin{pmatrix}
B^{(1)}\int \left[\frac{\partial f_\theta(x)}{\partial \theta_1} \otimes I_r\right] 
D^{(2)}(g(x),f_\theta(x)) d\mu(x) &  O \\ 
   	   O & O 
     	   \end{pmatrix},
$$
Now, note that $rank(H(\widetilde{\theta^g })) = r-1 $ and so we can not apply Theorem \ref{THM:IF_RMDE} 
directly to obtain the influence function of the RMDE in this case. However, we can restart with the 
set of equations  (\ref{EQ:IF_eq1}) and (\ref{EQ:IF_eq2}) with 
$\theta=\widetilde{\theta^g }= (\phi(\widetilde{\beta^g}), ~ \widetilde{\theta_2^g })^T$ and 
then solve those equations for the $IF$. For, let us partition the influence function 
$IF(y, \widetilde{T}_\rho, G) $ of $\widetilde{T}_\rho$ in terms of that of the functionals 
$\widetilde{T}_{\rho,1}$ and $\widetilde{T}_{\rho,2}$ corresponding to $\theta_1$ and $\theta_2$
 respectively as 
 $$
 IF(y, \widetilde{T}_\rho, G) = \begin{pmatrix}
      	    IF(y, \widetilde{T}_{\rho,1}, G) \\ 
      	   IF(y, \widetilde{T}_{\rho,2}, G)
      	    \end{pmatrix}.
 $$
 Now using the special form of the matrices for this case, those four equations simplifies to
\begin{eqnarray}
&& [BN_{11}(\widetilde{\theta^g})B^T]IF(y, \widetilde{T}_{\rho,1}, G) 
+ [BN_{12}(\widetilde{\theta^g })]IF(y, \widetilde{T}_{\rho,2}, G)   
\nonumber \\&&~~~
+ B^{(1)}\left\{\int \left[\frac{\partial f_\theta(x)}{\partial \theta_1} \otimes I_r\right] 
D^{(2)}(g(x),f_\theta(x)) d\mu(x) \right\} IF(y, \widetilde{T}_{\rho,1}, G)  
= B\left[\xi_1(\widetilde{\theta^g }) -  M_1(y;\widetilde{\theta^g })\right] , \label{EQ:eq1_RMDE2} 
\end{eqnarray}
\begin{eqnarray}
[N_{21}(\widetilde{\theta^g })B^T]IF(y, \widetilde{T}_{\rho,1}, G) + N_{22}(\widetilde{\theta^g})IF(y, \widetilde{T}_{\rho,2}, G)  = \left[\xi_2(\widetilde{\theta^g }) -  M_2(y;\widetilde{\theta^g })\right], \label{EQ:eq2_RMDE2} 
\end{eqnarray}
\begin{eqnarray}
[I_r - B^T] IF(y, \widetilde{T}_{\rho,1}, G) = 0. \label{EQ:eq3_RMDE2} 
\end{eqnarray}
Now from Equation (\ref{EQ:eq2_RMDE2}), we get 
\begin{equation}
IF(y, \widetilde{T}_{\rho,2}, G)  = N_{22}(\widetilde{\theta^g })^{-1}\left[\xi_2(\widetilde{\theta^g }) -  M_2(y;\widetilde{\theta^g })\right] -  N_{22}(\widetilde{\theta^g })^{-1}N_{21}(\widetilde{\theta^g })IF(y, \widetilde{T}_{\rho,1}, G). \label{EQ:IF2_RMDE2}
\end{equation}
Using this, Equation (\ref{EQ:eq1_RMDE2}) further simplifies to
\begin{eqnarray}
&& \left\{B\left[N_{11}(\widetilde{\theta^g }) 
- N_{12}N_{22}(\widetilde{\theta^g })^{-1}N_{21}(\widetilde{\theta^g })\right] 
+ B^{(1)}\int \left[\frac{\partial f_\theta(x)}{\partial \theta_1} \otimes I_r\right] 
D^{(2)}(g(x),f_\theta(x)) d\mu(x)\right\}IF(y, \widetilde{T}_{\rho,1}, G)  \nonumber \\
&& ~~~~~ = B\left\{\left[\xi_1(\widetilde{\theta^g }) -  M_1(y;\widetilde{\theta^g })\right] 
-  N_{12}(\widetilde{\theta^g })N_{22}(\widetilde{\theta^g })^{-1}\left[\xi_2(\widetilde{\theta^g }) 
-  M_2(y;\widetilde{\theta^g })\right]\right\}. \label{EQ:IF1_RMDE2}
\end{eqnarray}
We need to solve above for the first partition $IF(y, \widetilde{T}_{\rho,1}, G) $ subject to 
$B^T IF(y, \widetilde{T}_{\rho,1}, G) = IF(y, \widetilde{T}_{\rho,1}, G)$ and then use 
Equation \ref{EQ:IF2_RMDE2} to get the remaining second partition $IF(y, \widetilde{T}_{\rho,2}, G)$ of the IF.

In particular, if we have $N_{12}(\theta)=O$, then the estimators $\widetilde{\theta_1^g }$ and 
$\widetilde{\theta_2^g }$ becomes asymptotically independent and their influence functions also become 
independent of each other. The influence function of $\widetilde{\theta_2^g }$ becomes 
\begin{equation}
IF(y, \widetilde{T}_{\rho,2}, G)  = N_{22}(\widetilde{\theta^g })^{-1}\left[\xi_2(\widetilde{\theta^g }) -  M_2(y;\widetilde{\theta^g })\right]. \label{EQ:IF2_RMDE20}
\end{equation}
It is easy to see that this is indeed of the same form as the corresponding influence function in the unrestricted case. And , the influence function of $\widetilde{\theta_1^g }$, in this case, is given by the solution of 
\begin{eqnarray}
&&\left\{BN_{11}(\widetilde{\theta^g }) +  B^{(1)}\int \left[\frac{\partial f_\theta(x)}{\partial \theta_1} 
\otimes I_r\right] D^{(2)}(g(x),f_\theta(x)) d\mu(x)\right\}IF(y, \widetilde{T}_{\rho,1}, G)  
\nonumber \\ && ~~~~~~~~~~~~~~~~~~~~~~~~~~~~~~~~~~~~~~~~~~~~~~~ 
= B\left\{\left[\xi_1(\widetilde{\theta^g }) -  M_1(y;\widetilde{\theta^g })\right]\right\},
~~~~~ \label{EQ:IF1_RMDE20}
\end{eqnarray}
subject to the restriction $B^T IF(y, \widetilde{T}_{\rho,1}, G) = IF(y, \widetilde{T}_{\rho,1}, G)$.

Now let us try to derive the influence function for our motivating case in this example; namely, the
p-variate normal model with mean $\mu$ and variance $\sigma^2 I_p$ with the restriction $\mu = \beta \mu_0$. 
Thus, here, $\phi_i(\beta)=\beta (\mu_0)_i$ for all $i=1,\cdots, p$. Hence we have $b_{ij} = constant$ for all $i, j$ and so $B^{(1)}=O$. Further, considering $\theta_1=\mu$ and $\theta_2=\sigma^2$,
in this case we have $N_{12}(\theta)=0$. Thus, from above, the influence function of $\widetilde{{\sigma^2}^g }$ is given by 
 $$IF(y, \widetilde{{\sigma^2}^g }, G)  = N_{22}(\widetilde{\theta^g })^{-1}\left[\xi_2(\widetilde{\theta^g }) -  M_2(y;\widetilde{\theta^g })\right].$$ 
And the influence function of  $\widetilde{{\mu}^g }$ is then a solution of 
\begin{equation}
B N_{11}(\widetilde{\theta^g })IF(y, \widetilde{{\mu}^g }, G)  = B\left[\xi_1(\widetilde{\theta^g }) -  M_1(y;\widetilde{\theta^g })\right], \label{EQ:IF1_RMDE20_norm}
\end{equation}
subject to the restriction $B^T IF(y, \widetilde{{\mu}^g }, G) = IF(y, \widetilde{{\mu}^g }, G)$.
Thus we will get a non-zero influence function of $\widetilde{{\mu}^g }$ if the matrix $B^T$ 
has one of its eigenvalue as $1$ and in that case the influence function is given by that eigenvalue 
of $B^T$ corresponding to the eigenvalue $1$ which satisfies the equation (\ref{EQ:IF1_RMDE20_norm}). 
After simplification, in that case, the influence function must be of the form 
$$ IF(y, \widetilde{{\mu}^g }, G)  = B^T N_{11}(\widetilde{\theta^g })^{-1}\left[\xi_1(\widetilde{\theta^g }) -  M_1(y;\widetilde{\theta^g }) + v\right], $$
 where $v$ is a vector in the null-space of the matrix $B$.
 
 	For the special choice $\mu_0 =(1, \cdots, 1)^T$, we have $b_{ij}=1$ for all $i, j$ so that the matrix $B$ 
 	does not have eigenvalue $1$ and hence $IF(y, \widetilde{{\mu}^g }, G) = 0$.   
 \hfill{$\square$}

\section{Applications : Some Particular Divergences}\label{SEC:part_div}

Based on the general results obtained in the two previous sections, one can describe the influence function 
analysis and the asymptotic distributions of any MDE or RMDE provided, she can prove only their $\sqrt n$-consistency.
In this section, we will apply those results for some common divergence measures and common model family. 
Throughout this section, we will assume some common notations from the likelihood theory as, 
$L(\theta; \Theta) = \ln f_\theta(x)$ for all $\theta \in \Theta$ is the likelihood function, 
$u_\theta(x) = \nabla L(\theta; \Theta)$ is the  the likelihood score function,
$I(\theta)= E_{f_\theta}[u_\theta(X) u_\theta(X)^T]$ is the fisher information matrix. 
Also, we will define similar quantities under a proper subspace $\Theta_0 \subset \Theta$ (different by the
 restrictions $h(\theta) =0$ as, $L(\theta; \Theta_0)$ being the restriction of $L(\theta; \Theta)$ onto the subspace
 $\Theta_0$, $u_\theta^0(x)=\nabla L(\theta; \Theta_0)$ and $I^0(\theta)=E_{f_\theta}[u_\theta^0(X) u_\theta^0(X)^T]$.

\subsection{Disparity Measures}

One of the most popular family of divergences is the disparity family (Lindsay, 1994) that yields fully efficient and 
robust estimators upon minimization. It is defined in terms of a non-negative thrice differentiable strictly convex function $\phi$ on $[-1, \infty)$ with $\phi(0)=0$ and $\phi'(0)=0$, called the disparity generating function, 
as 
$$
\rho(g, f) = \int \phi(\delta)f,
$$
with $\delta = g/f - 1$. It is of the form of general divergences defined in Equation (\ref{EQ:part_div}) with 
$D(a, b) = \phi\left(\frac{a}{b} - 1\right)$ so that we can apply all the results derived above. 
Using the same notations, we have, 
$$
M(y; \theta) = - A'(\delta)u_\theta(y),
$$
and
$$
N(\theta) = \int \left[A'(\delta)u_\theta u_\theta^T g - A(\delta)\nabla_2 f_\theta \right],
$$
where, the function $A(\delta)$, defined as $A(\delta) = C'(\delta)(\delta + 1) - C(\delta)$, is known as the 
Residual Adjustment Function in the context of minimum disparity estimation and plays a crucial role in 
its robustness (Lindsay, 1994). Thus, using (\ref{EQ:IF_MDE}), the influence function of the 
minimum disparity estimator is given by 
\begin{eqnarray}
 IF(y, T_\rho, G)  = N(\theta^g)^{-1} [A'(\delta)u_\theta(y) - E_g[A'(\delta)u_\theta(y)]],
\end{eqnarray}
which is the same as obtained by Lindsay (1994) independently. In particular, at the model $g=f_{\theta_0}$,
this influence function simplifies to $I(\theta_0)^{-1}u_{\theta_0}(y)$ which is independent of the disparity 
generating function $C(\cdot,\cdot)$ and so is same as that of the MLE. This is unbounded function for most of 
the common model families.

Now, let us consider the restricted minimum disparity estimation under the restrictions $h(\theta) = 0$.
Using the notations of Section \ref{SEC:IF_RMDE_general}, 	it is easy to see that 
$
M_0(y; \theta) = - A'(\delta)u_\theta^0(y),
$
and
$$
N_0(\theta) = \int \left[A'(\delta)u_\theta^0 (u_\theta^0)^T g 
- A(\delta)\nabla_2[ f_\theta]\big|_{\Theta_0} \right].
$$
Then, we can derive the influence function of the restricted minimum disparity estimators from Theorem \ref{THM:IF_RMDE}
and above simplified expressions. However, the interesting case is when true density belongs to the model family, i.e.
$g=f_{\theta_0}$. In that case, we will have $M_0(y; \theta_0) = - u_\theta^0(y)$, $\xi_0(\theta_0)=0$ and
$N_0(\theta_0)=I^0(\theta_0)$. Then, we get the simple expression of the restricted minimum disparity estimator 
$\widetilde{T}_C$ corresponding to the disparity generated by $C(\cdot,\cdot)$ as
\begin{eqnarray}
IF(y, \widetilde{T}_C, F_{\theta_0})  = \left[I^0({\theta_0 })^2  + 
   H({\theta_0 }) H({\theta_0 })^T\right]^{-1} I^0({\theta_0 }) u_\theta^0(y) .
\end{eqnarray} 
Note that the above expression is independent of the choice of the disparity generating function and hence it also 
gives the influence function of the Restricted Maximum Likelihood Estimators. 

Further, it will help us to derive asymptotic distribution of the restricted minimum disparity estimators 
$\widetilde{\theta}_C$ including that of the restricted maximum likelihood estimators. Following the argument of
 Lindsay (1994), one can easily prove the $\sqrt{n}$-consistency of the restricted minimum divergence estimators. 
 Then, as pointed out in Remark \ref{REM:asymp_RMDE}, the asymptotic distribution of 
 $\sqrt{n}(\widetilde{\theta}_C - \theta_0)$, at the model $g=f_{\theta_0}$, is normal with mean zero and 
 variance given by
 $$
 \left[I^0({\theta_0 })^2  + 
    H({\theta_0 }) H({\theta_0 })^T\right]^{-1} [I^0({\theta_0 })]^3\left[I^0({\theta_0 })^2  + 
       H({\theta_0 }) H({\theta_0 })^T\right]^{-1}.
 $$
This expression coincides with the asymptotic distribution of restricted maximum likelihood estimators obtained 
independently from the likelihood theory. 
Hence it provides a justification of our general results obtained in this paper.

\subsection{Density Power Divergence}

In the recent decades, arguably the most popular divergence measure in the context of the robust minimum divergence 
estimation is the Density Power Divergence (Basu et al., 1998). The increasing popularity of this divergence is 
mainly due to the fact that corresponding minimum divergence estimation does not require any kernel smoothing for the 
continuous models; which is a major drawback of disparity measures. 
The density power divergence is defined in terms of a non-negative tuning parameter $\alpha$ as 
$$
	\rho_\alpha(g,f) = \int  f^{1+\alpha} - \frac{1+\alpha}{\alpha} \int f^\alpha g + 
	\frac{1}{\alpha} \int g^{1+\alpha}  ~~~~~~~~~~~ \mbox{ for } \alpha > 0, 
$$ and 
$$
\rho_0(g,f) = \lim_{\alpha \rightarrow 0} d_\alpha(g,f) = \int g \log(g/f) .
$$
Note that the case of $\alpha=0$ gives the likelihood disparity and so the influence function of the corresponding 
minimum divergence estimator is already discussed in previous subsection. Let us now consider the case $\alpha>0$.
Interestingly, for any given fixed $\alpha>0$, this divergence also belongs to the general family of divergence
 defined in Equation (\ref{EQ:part_div}) with 
$$
D(a,b) = b^{1+\alpha} - \left(\frac{1+\alpha}{\alpha}\right)b^\alpha a + \frac{1}{\alpha} a^{1+\alpha}.
$$
Now, we can apply all the results derived above for the density power divergences where, 
$
M(y; \theta) = - (1+\alpha)u_\theta(y)f_\theta^\alpha(y),
$
and
$$
N(\theta) = (1+\alpha) \int \left[u_\theta u_\theta^T f_\theta^{1+\alpha} + (i_\theta - \alpha u_\theta u_\theta^T)(g - f_\theta)f_\theta^\alpha \right] = (1+\alpha)J_\alpha(\theta),
$$
where, $i_\theta = - \nabla u_\theta$.
Thus, from (\ref{EQ:IF_MDE}), the influence function of the 
minimum density power divergence estimator $T_\alpha$ can be written as 
\begin{eqnarray}
 IF(y, T_\alpha, G)  = J_\alpha(\theta^g)^{-1} [u_\theta(y)f_\theta^\alpha(y) - E_g[u_\theta(X)f_\theta^\alpha(X)]],
\end{eqnarray}
which exactly as derived in Basu et al.~(1998). In particular, if we assume $g=f_{\theta_0}$,
then the influence function becomes 
\begin{equation}\label{EQ:IF_DPD0}
 IF(y, T_\alpha, F_{\theta_0})  = \left(\int u_{\theta_0} u_{\theta_0}^T f_{\theta_0}^{1+\alpha}\right)^{-1} 
 \left[u_{\theta_0}(y)f_{\theta_0}^\alpha(y) - \int u_{\theta_0} f_{\theta_0}^{1+\alpha}\right].
\end{equation}
This influence function is bounded for all $\alpha>0$ and most of the common model families.

Next, we will consider the restricted minimum density power divergence estimation under the 
restrictions $h(\theta) = 0$. Again, we use the notations of Section \ref{SEC:IF_RMDE_general} so that
$
M_0(y; \theta) = - (1+\alpha)u_\theta^0(y)f_\theta^\alpha(y),
$
and
$$
N_0(\theta) = (1+\alpha) \int \left[u_\theta^0 (u_\theta^0)^T f_\theta^{1+\alpha} 
+ (i_\theta^0 - \alpha u_\theta^0 (u_\theta^0)^T)(g - f_\theta)f_\theta^\alpha \right] 
= (1+\alpha)J_\alpha^0(\theta),
$$
with $i_\theta^0 = - \nabla u_\theta^0$.
Then, Theorem \ref{THM:IF_RMDE} gives us the expression of the influence function of the restricted 
density power divergence estimators. In particular, if $g=f_{\theta_0}$, then the influence function of the 
restricted minimum disparity estimator $\widetilde{T}_\alpha$ simplifies to 
\begin{eqnarray}\label{EQ:IF_RMDPDE0}
IF(y, \widetilde{T}_\alpha, F_{\theta_0})  = \Psi(\theta_0)^{-1} 
\left(\int u_{\theta_0}^0 (u_{\theta_0}^0)^T f_{\theta_0}^{1+\alpha}\right) u_{\theta_0}^0(y) f_{\theta_0}^\alpha(y),
\end{eqnarray} 
where $\Psi(\theta) = \left[\left(\int u_\theta^0 (u_\theta^0)^T f_\theta^{1+\alpha}\right)^2  + 
\frac{1}{(1+\alpha)^2} H({\theta }) H({\theta})^T\right]$; 
again this influence function is generally bounded for all $\alpha>0$.

Finally, we can derive the asymptotic distribution of the restricted minimum density power divergence estimators 
$\widetilde{\theta}_\alpha$ from Remark \ref{REM:asymp_RMDE}. The $\sqrt{n}$-consistency of the 
restricted minimum density power divergence estimators follows from a modification of the argument of Basu et al.~
(1998) used to prove the same for minimum density power divergence estimators. So, if we have $g=f_{\theta_0}$, then 
 asymptotic distribution of  $\sqrt{n}(\widetilde{\theta}_\alpha - \theta_0)$ is normal with mean zero and 
 asymptotic variance 
\begin{eqnarray}
 \Psi(\theta_0)^{-1} \left(\int u_{\theta_0}^0 (u_{\theta_0}^0)^T f_{\theta_0}^{1+\alpha}\right) 
 Var_{f_{\theta_0}}[u_{\theta_0}^0(Y) f_{\theta_0}^\alpha(Y)]\left(\int u_{\theta_0}^0 (u_{\theta_0}^0)^T f_\theta^{1+\alpha}\right) 
 \Psi(\theta_0)^{-1}. \nonumber
\end{eqnarray}

\subsection{$S$-Divergence}

We will now consider a recent family of divergences, namely the $S$-Divergence Family, developed by
 Ghosh et al.~(2013). This is a general super-family containing both the density power divergence 
 (Basu et al., 1998) and the Cressie-Read family of power divergences (Cressie and Read, 1984) and also contains many other useful divergences. It is defined in terms of tqo parameters 
 $\lambda \in \mathbb{R}$ and $\alpha \geq 0$ as 
\begin{equation}
\rho(g,f) = S_{(\alpha, \lambda)}(g,f) =  \frac{1}{A} ~ \int ~ f^{1+\alpha}  -   \frac{1+\alpha}{A B} ~ 
\int ~~ f^{B} g^{A}  + \frac{1}{B} ~ \int ~~ g^{1+\alpha}, ~~ \alpha \in [0, 1], ~~\lambda \in {\mathbb R},
\label{EQ:S_div_gen}
\end{equation}
where, $A = 1+\lambda (1-\alpha)$ and  $B = \alpha - \lambda (1-\alpha)$. For either $A=0$ or $B=0$, 
it is defined by the corresponding continuous limit of divergences [See Ghosh et al.~(2013) for details]. Again, this large family of divergence can be written in the form of equation
 (\ref{EQ:part_div}) with 
$$
D(a, b) = \frac{1}{A} b^{1+\alpha} - \left(\frac{1+\alpha}{\alpha}\right)b^B a^A + \frac{1}{B} a^{1+\alpha}.
$$
Then, we have $M_0(y; \theta) = - (1+\alpha)u_\theta(y)f_\theta^B(y) g^A(y)$,
and
$$
N(\theta) = \frac{(1+\alpha)}{A} \int \left[Au_\theta u_\theta^T f_\theta^{1+\alpha} + (i_\theta^0 - B u_\theta u_\theta^T)(g^A - f_\theta^A)f_\theta^B \right] = (1+\alpha) J_{(\alpha, \lambda)}.
$$
Then, we ge the influence function of the minimum $S$-divergence estimator from equation \ref{EQ:IF_MDE} as given by 
\begin{eqnarray}
 IF(y, T_{(\alpha, \lambda)}, G)  = J_{(\alpha, \lambda)}(\theta^g)^{-1} [u_\theta(y)f_\theta^B(y) g^A(y) - E_g[u_\theta(X)f_\theta^B(X) g^A(X)]],
\end{eqnarray}
which is again the exactly same as obtained in Ghosh et al.~(2013). For the special case $g=f_{\theta_0}$,
this influence function coincides with that of the density power divergence given by equation (\ref{EQ:IF_DPD0}).

Finally, the influence function of the restricted minimum $S$-divergence estimators 
under the restrictions $h(\theta)=0$ can be derived from Theorem \ref{THM:IF_RMDE}. 
It is then easy to see that, at the model $g=f_{\theta_0}$, 
 the influence function of the restricted minimum $S$-divergence estimators coincides with that of the restricted 
 density power divergence estimators derived in equation (\ref{EQ:IF_RMDPDE0}).

\section{Conclusion}
This work present the derivation of the influence function of the restricted and unrestricted 
minimum divergence estimators for a general class of density based divergences.
It will help researchers to derive the robustness properties of any minimum divergence estimators 
under several restrictions on the parameters. As an example, we have examined the same
for some popular minimum divergence estimators, namely the  disparity, density power divergence and 
$S$-divergence family; we have also presented an example with a set of linearly dependent restrictions
for general model family.
Further, this paper gives us several directions for future works including the influence function of 
more general class of divergences that are possibly  based on the distribution functions; 
author want to solve the related problems in subsequent researches.
\\\\

\noindent \textbf{Acknowledgment:}	This work is a part of the author's Ph.D. dissertation under 
Professor Ayanendranath Basu. Author want to thank Prof. Basu for his sincere guidance and comments
regarding this paper.


\begin{thebibliography}{}

\bibitem{al66}
Ali, S. M. and Silvey, S. D. (1966). A general class of coefficients of divergence of 
one distribution from another. {\it J. Roy. Statist. Soc.}, {\bf B 28}, 131--142.

\bibitem{bhhj98}
Basu, A., Harris, I. R., Hjort, N. L. and Jones, M. C. (1998). Robust and efficient 
estimation by minimizing a density power divergence. {\it Biometrika}, {\bf 85}, 549--559.

\bibitem{bsp11}
Basu, A., Shioya, H. and Park, C. (2011). {\it Statistical Inference: The Minimum Distance Approach}. Chapman $\&$ Hall/CRC.

\bibitem{bmmp13}
Basu, A., Mandal, A., Martin, N. and Pardo, L. (2013). Restricted Power Divergence Tests for Composite Null Hypotheses {\it Technical Report, Bayesian and Interdisciplinary Research Unit, Indian Statistical Institute, Kolkata, India}

%
\bibitem{c63}
Csisz\'{a}r, I. (1963). Eine informations theoretische Ungleichung und ihre Anwendung auf den 
Beweis der Ergodizitat von Markoffschen Ketten. {\em Publ. Math. Inst. Hungar. Acad. Sci.}, {\bf 3}, 85--107.
%
\bibitem{c67a}
Csisz\'{a}r, I. (1967a). Information-type measures of difference of probability distributions and indirect observations. {\em Studia Scientiarum Mathematicarum Hungarica}, {\bf 2}, 299--318.

\bibitem{c67b}
Csisz\'{a}r, I. (1967b). On topological properties of f-divergences. {\em Studia Scientiarum Mathematicarum Hungarica}, {\bf 2}, 329--339.

\bibitem{cr84}
Cressie, N. and T. R. C. Read (1984). Multinomial goodness-of-fit tests. {\it J. Roy. Statist. Soc.},
{\bf B 46}, 440--464.


\bibitem{gmb13}
Ghosh, A., Maji, A., Basu, A. (2013).  Robust Inference based on Divergences in Reliability Systems. {\em Applied Reliability Engineering and Risk Analysis: Probabilistic Models and Statistical Inference}. Eds, Inference, Ilia Frenkel, Alex, Karagrigoriou, Anatoly Lisnianski and Andre Kleyner,
Dedicated to the Centennial of the birth of Boris Gnedenko. 
John Wiley $\&$ Sons, Ltd . 290--307


\bibitem{h68}
\textsc{Hampel, F. R.} (1968). \textit{Contributions to the Theory of Robust Estimation}. Ph.D. Thesis, University of California, Berkeley, U.S.A.

\bibitem{h74}
\textsc{Hampel, F. R.} (1974). The influence curve and its role in robust estimation. \textit{J. Amer. Statist Assoc.}, {\bf 69}, 383--393.

\bibitem{hrrs86}
\textsc{Hampel, F. R.}, \textsc{Ronchetti, E. M.}, \textsc{Rousseeuw, P. J.} and \textsc{Stahel W. A.} (1986). \textit{Robust Statistics: The Approach Based on Influence Functions}. John Wiley $\&$ Sons, New York. 


\bibitem{l94}
Lindsay, B. G. (1994). Efficiency versus robustness: The case for minimum 
Hellinger distance and related methods. {\em Ann. Statist.}, {\bf 22}, 1081--1114.

\bibitem{p06}
Pardo, L. (2006). {\em Statistical Inference based on Divergences}. CRC/Chapman-Hall. 

\bibitem{pmbp13}
Patra, S., Maji, A., Basu, A., Pardo, L. (2013).The Power Divergence and the 
Density Power Divergence Families : the Mathematical Connection. {\it Sankhya B}, {\bf 75}, 16--28.

\bibitem{s87}
\textsc{Simpson, D. G.} (1987). Minimum Hellinger distance estimation for the analysis of count data. \textit{J. Amer. Statist. Assoc.} {\bf 82}, 802--807.

\bibitem{vi89}
Vajda, I. (1989). {\it Theory of Statistical Inference and Information}. Dordrecht: Kluwer Academic. 

\end{thebibliography}
\end{document}